\documentclass[11pt]{article}

\usepackage{amssymb}
\usepackage{amsmath}
\usepackage{epsfig}
\usepackage{color}

\newtheorem{theorem}{Theorem}[section]
\newtheorem{lemma}[theorem]{Lemma}

\newtheorem{example}[theorem]{Example}

\renewcommand{\baselinestretch}{1.5}

\title{Orthogonal cycle systems with cycle length less than $10$}

\author{{\bf{Selda K\"{u}\c{c}\"{u}k\c{c}if\c{c}i} }\\
   email: {\tt skucukcifci@ku.edu.tr}\\
{\bf{Emine \c{S}ule Yaz\i c\i} }\\
    email: {\tt eyazici@ku.edu.tr}    \\
    Department of   Mathematics, Ko\c{c} University, Istanbul,  Turkey}

\date{ }

\begin{document}
\maketitle\thispagestyle{empty}

\def\baselinestretch{1.5}\small\normalsize
\begin{abstract}

\end{abstract}

An $H$-decomposition of  $G$ is a partition of the edge-set of $G$ into subsets, where each subset induces a copy of the graph $H$. A $k$-orthogonal $H$-decomposition of a graph $G$ is a set of $k$ $H$-decompositions of $G$, such that any two copies of $H$ in distinct $H$-decompositions intersect in at most one edge. When $G=K_v$ we call the $H$-decomposition an $H$-system of order $v$. In this paper we consider the case $H$ is an $l$-cycle and construct a pair of orthogonal $l$-cycle systems for all admissible orders when $l=5,6,7, 8\ or\ 9$, except $(l,v)=(7,7)$ and $(l,v)=(9,9)$.

\noindent {\small AMS classification: $05C38$, $05B15$, $05B30$.\\
 Keywords: Orthogonal cycle decompositions; Heffter arrays}

\section{Introduction}

Orthogonal combinatorial structures especially; orthogonal Latin squares has a history dating back to the famous Euler's 36 officer problem. Please see for example the vast literature on orthogonal Latin squares in \cite{Handbook} and the survey of orthogonal factorizations in \cite{AHL}. We will consider orthogonal cycle decompositions in this paper.

An $H$-decomposition of  $G$ is a partition of the edge-set of $G$ into subsets, where each subset induces a copy of the graph $H$. A $k$-orthogonal $H$-decomposition of a graph $G$ is a set of $k$ $H$-decompositions of $G$, such that any two copies of $H$ in distinct $H$-decompositions intersect in at most one edge. When $G=K_v$ we call the $H$-decomposition an $H$-system of order $v$. Copies of $H$ are called blocks.

In \cite{CY}, it was shown asymptotically that for sufficiently large $v$, $k$-orthogonal $K_r$-decompositions of $K_v$ exists whenever a $K_r$-decomposition of $K_n$ exists. This result was generalized in \cite{CY2} to the general graph $H$ which we will refer as an $H$-system of order $v$ . These structures are closely related to super-simple and completely-reducible super-simple designs. See \cite{Handbook} for definitions and related results.

It was established in \cite{AG,Bur, sajna} that when $H$ is an $l$-cycle then a decomposition of order $v$ exists if and only if $v$ is odd, $3\leq l\leq v$ and $v(v-1)\equiv 0\ (mod \ 2l)$. Such a decomposition is called an $l$-cycle system of order $v$. When $H$ is a $3$-cycle,  two $3$-cycle systems are orhogonal if and only if the cycle systems are disjoint. Therefore, there are at most $v-2$ mutually orthogonal $3$-cycle systems of order $v$. A set of systems which meets this bound forms a large set of disjoint Steiner triple systems and they exist for all admissible $v \neq 7$ \cite{Lu1,Lu2,Teirlinck}.

\begin{example} \label{exK7}
A pair of orthogonal Steiner triple systems of order $7$ exists.
\end{example}
Let Let $V=\mathbb{Z}_{7}$, $C_1$ and $C_2$, each, be a collection of triples obtained by developing the base triples $(1,2,4)$ and $(1,3,4)$, respectively modulo $7$.  Then $(V,C_1)$ and $(V,C_2)$ form a pair of orthogonal Steiner triple systems of order $7$. \hfill $\square$

In \cite{BCP}, authors focus on the case $v\equiv 1 \ (mod  \ 2l)$ when $l$ is even and proved the following result.

\begin{theorem}\cite{BCP}
Suppose $\mu(v,l)$ represent the maximum possible number of mutually orthogonal $l$ cycle systems of order $v$. Let $l\geq 4$  is even, $v\equiv\ 1\ (mod\ 2l)$ and $N=\frac{v - 1}{2l}$. Then $\mu(v,l)\geq \frac{N}{al+b}-1$, where $(a, b)=(4,-2)$ if  $ l\equiv\ 0\ (mod\ 4)$ and; $(a,b)=(24,-18)$ if $ l\equiv\ 2\ (mod\ 4)$.
\end{theorem}

Also in the same paper it was noted that $\mu(v,l)\leq v-2$ in general and $\mu(v,l)\leq 1$ when $2l^2>v(v-1)$.

On the other hand, orthogonal cycle systems naturally arise from face 2-colourable embeddings of graphs on surfaces, where each pair of faces share at most one edge and each edge belongs to exactly two faces. So first results on determining the spectrum of 2-orthogonal cycle systems make use of Heffter arrays that are constructed for the purpose of obtaining these type of  2-colourable embeddings of the complete graph on high genius surfaces.
 A {\em Heffter array} $H(m,n;s,l)$ is an $m\times n$ array of integers, such that:
\begin{itemize}
\item each row contains $s$ filled cells and each column contains $l$ filled cells;
\item the elements in every row and column sum to $0 \ (mod \ 2ml+1)$; and
\item for each integer $1\leq x\leq ms$, either $x$ or $-x$ appears in the array.
\end{itemize}

A Heffter array is {\em square}, if $m = n$ and necessarily $s = l$, and is denoted by $H(m;l)$. A square Heffter array $H(m;l)$ is said to be a simple array if, for each row and
column, the entries may be cyclically ordered so that all partial sums, $$\alpha_i=\displaystyle\sum_{j=0}^{i} {a_j} \ mod \ (2ml+1)$$   where $a_0,...,a_{l-1}$ are the elements in that row or column, are distinct modulo
$2ml + 1$ for $i=0,\dots, l-1$. One may easily check that any Heffter array with $l\leq 5$ is simple. The results in \cite{Archdeacon} implies the following theorem as part of a more general result; but it was stated explicitly in \cite{CMPP} including the case for the cocktail party graphs.

\begin{theorem} If $H(m; l)$ is a Heffter array with a simple ordering, then there exists a pair of orthogonal cyclic $l$-cycle systems of order $2ml+1$.
\end{theorem}

The $H(m;6)$'s constructed in \cite{CDDY} are simple. Thus the following is implied by the existing literature on Heffter arrays.

\begin{theorem} \label{thHeffter} \cite{ADDY, BCDY, CDDY, CMPP, DW} Let $ m \geq l \geq 3$ be positive integers. Then  there exists a pair of orthogonal $l$-cycle systems of order $2ml+1$ when
\begin{itemize}
 \item  $3\leq l \leq 6$;
 \item $l \in \{7, 9\}$ and $ml\equiv 3\ (mod\ 4)$;
 \item $l\equiv 0 \ (mod\ 4)$;
 \item $m\equiv 1\ \ (mod\ 4)$ and $l\equiv 3 \ (mod\ 4)$;
 \item $m\equiv 0\ (mod\ 4)$ and $l\equiv 3 \ (mod\ 4)$ $($for large enough $m)$.
\end{itemize}
\end{theorem}

If they satisfy an extra condition on the orderings of the entries of the Heffter array, these orthogonal cycle systems in turn biembed to yield a face 2-colourable embedding on an orientable surface. (See for example \cite{CDY,DM}). For other existence results and generalizations of Heffter arrays, see \cite{ABD, CDDY, CDP2, ItalE2, CP, ItalE3, MP1, MP2, MP3}. 

All the orthogonal $l$-cycle systems obtained using Heffter arrays given in Theorem \ref{thHeffter} and the cyclic systems constructed in \cite{BCP} are of order $v=2ml+1$ for some integer $m$, hence only deals with the case $v\equiv 1\ (mod \ 2l)$ for $m\geq l$. In this paper  we will show that a pair of orthogonal $l$-cycle systems of order $n$ exists whenever an $l$-cycle system of order $n$ exists for the cycle lenghts $3\leq l\leq 9$.

When $l=3$, as stated earlier a large set of Steiner triple systems exists for all $v\neq 7$ and the unique $STS(7)$ has an orthogonal pair hence a pair of orthogonal $3$-cycle systems exist for all $v\equiv 1,3 \ (mod \ 6)$.  When $l=4$, a $4$-cycle system exists if and only if $v\equiv 1 \ (mod\ 8)$. A pair of orthogonal $4$-cycle systems is obtained in \cite{BCP} for each of these admissible orders. We will construct a pair of orthogonal $l$-cycle systems for all admissible orders when $l=5,6,7, 8\ or\ 9$, except $(l,v)=(7,7)$ and $(l,v)=(9,9)$.
.

The followings are the definitions of the necessary structures for our constructions.

A $k-$GDD is a triple $(V,\cal G, B)$, where $V$ is a finite set of vertices, $\cal G$ $=\{G_1,G_2,\ldots,G_n\}$ is a partition of $V$ into subsets, called {\em groups} and $\cal B$ is a collection of
isomorphic copies of $K_k$, which partitions the edges of $K_{g_1, g_2,\ldots, g_n}$ on $V$, where $|G_i| = g_i$. If for $i=1,2,\ldots, t$, there are $u_i$ groups of size $g_i$, we
say that the $k-$GDD is of type $g_1^{u_1}g_2^{u_2}\ldots g_t^{u_t}$.

We will denote all cycles with length $l$ and edges $\{a_1,a_2\},\{a_2,a_3\},...,\{a_{l-1},a_l\},\{a_l,a_1\}$ by $(a_1,a_2,...,a_l)$.

An automorphism of an $STS(v)$ $(V, T)$ is a bijection $\alpha : V \rightarrow V$ such that $t=\{x,y,z\} \in T$ if and only if $t\alpha = \{x\alpha ,y\alpha ,z\alpha \}\in T$. An $STS(v)$ is cyclic if it has an automorphism that is a permutation consisting of a single cycle of length $v$.

\section {Orthogonal $5$-cycle systems}

The spectrum of $5$-cycle systems; that is the set of all $v$ such that a $5$-cycle system of order $v$ exists, is the set of all $v\equiv 1$ or $5$ (mod $10$). There exists a pair of orthogonal cyclic $5$-cycle systems of order $10k+1$ for all $k\geq 5$ \cite{CMPP}. But we will give a construction that works for both of the orders $10k+1$ and $10k+5$ for all $k\geq 3$.

Since a pair of orthogonal $5$-cycle systems of order $11$ is given in \cite{BCP}, we will start by giving the remaining pairs of orthogonal $5$-cycle decompositions we will need in our constructions. All systems of small orders that are not used in our recursive constructions throughout the paper will be given in the Appendix.

\begin{example} \label{ex15}
A pair of orthogonal $5$-cycle systems of order $15$ exists.
\end{example}

Let $V=\mathbb{Z}_{15}$ and

$C_1= \{(0, 1, 2, 3, 4)$, $(0, 2, 12, 8, 5)$, $(0, 3, 1, 4, 6)$, $(0, 7, 1, 5, 9)$, $(0, 8, 1, 6, 10)$, $(0, 11, 1, 9, 12)$, $(0, 13, 1, 10, 14)$, $(1, 12, 3, 5, 14)$, $(2, 4, 5, 6, 7)$, $(2, 5, 7, 3, 6)$, $(2, 8, 3, 9, 10)$,
$(2, 9, 4, 7, 11)$, \linebreak
$(2, 13, 3, 11, 14)$, $(3, 10, 4, 8, 14)$, $(4, 11, 5, 10, 12)$, $(4, 13, 5, 12, 14)$, $(6, 8, 7, 9, 11)$, $(6, 9, 8, 13, 14)$, \linebreak
$(6, 12, 11, 10, 13)$, $(7, 10, 8, 11, 13)$, $(7, 12, 13, 9, 14)\}$,

$C_2= \{(0, 1, 3, 5, 2)$, $(0, 3, 2, 4, 7)$, $(0, 4, 1, 5, 8)$, $(0, 5, 4, 3, 6)$, $(0, 9, 1, 2, 10)$, $(0, 11, 2, 6, 13)$, \linebreak
$(0, 12, 1, 6, 14)$, $(1, 7, 2, 8, 10)$, $(1, 8, 3, 7, 11)$, $(1, 13, 2, 9, 14)$, $(2, 12, 3, 10, 14)$, $(3, 9, 4, 6, 11)$, \linebreak
$(3, 13, 5, 7, 14)$, $(4, 8, 6, 5, 12)$, $(4, 10, 5, 9, 13)$, $(4, 11, 9, 8, 14)$, $(5, 11, 8, 12, 14)$, $(6, 7, 8, 13, 10)$, \linebreak
$(6, 9, 7, 10, 12)$, $(7, 12, 11, 14, 13)$, $(9, 10, 11, 13, 12)\}$. Then $(V,C_1)$ and $(V,C_2)$ form a pair of orthogonal $5$-cycle systems of order $15$. \hfill $\square$

\begin{example} \label{ex15-5}
A pair of orthogonal $5$-cycle decompositions of $K_{15}\setminus K_5$ exists.
\end{example}

Let $V=\mathbb{Z}_{15}$ with the hole of size $5$ on $\{0,1,2,3,4\}$ and

$C_1= \{(0, 5, 11, 6, 10)$, $(0, 6, 5, 2, 14)$, $(0, 7, 1, 5, 8)$, $(0, 9, 1, 6, 12)$, $(0, 11, 1, 8, 13)$, $(1, 10, 2, 6, 13)$, \linebreak
$(1, 12, 2, 7, 14)$, $(2, 8, 3, 5, 9)$, $(2, 11, 3, 7 , 13)$, $(3, 6, 4, 5, 10)$, $(3, 9, 4, 7, 12)$, $(3, 13, 4, 8, 14)$, \linebreak
$(4, 10, 7, 5, 12)$, $(4, 11, 7, 6, 14)$, $(5, 13, 9, 12, 14)$, $(6, 8, 10, 11, 9)$, $(7, 8, 12, 10, 9)$, $(8, 9, 14, 13, 11)$, \linebreak
$(10, 13, 12, 11, 14)\}$,

$C_2= \{(0, 5, 1, 6, 13)$, $(0, 6, 2, 7, 8)$, $(0, 7, 3, 5, 10)$, $(0, 9, 2, 5, 11)$, $(0, 12, 1, 8, 14)$, $(1, 7, 4, 5, 9)$, \linebreak
$(1, 10, 3, 8, 11)$, $(1, 13, 3, 6, 14)$, $(2, 8, 4, 6, 10)$, $(2, 11, 4, 9, 12)$, $(2, 13, 4, 10, 14)$, $(3, 9, 6, 5, 14)$, \linebreak
$(3, 11, 6, 8, 12)$, $(4, 12, 7, 9 ,14)$, $(5, 7, 6, 12, 13)$, $(5, 8, 9, 11, 12)$, $(7, 10, 8, 13, 11)$, $(7, 13, 10, 12, 14)$, \linebreak
$(9, 10, 11, 14, 13))\}$. Then $(V,C_1)$ and $(V,C_2)$ form a pair of orthogonal $5$-cycle decompositions of $K_{15}\setminus K_5$. \hfill $\square$

\begin{lemma} \label{Ort5}
A pair of orthogonal $5$-cycle systems of order $v$ exists for all $v\equiv 1$ or $5$ $($mod $10)$.
\end{lemma}

\noindent {\bf Proof} \ The cases $v=11$ is given in \cite{BCP}, $v=15$ in Example \ref{ex15} and $v=21, 25$ in the Appendix.  The following construction works for all $v=10k+r$ and $k \geq 3$, when $r=1$ or $5$. Let $(Q, \circ )$ be a commutative quasigroup of order $2k$, for $k\geq 3$, with holes $H=\{h_1,h_2,...,h_k\}$, where $h_i=\{2i-1,2i\}$ for $i=1,2,...,k$. Let $V=\displaystyle \bigcup_{j=1}^r \{\infty_j\}\cup (Q\times \{1,2,3,4,5\})$, where $r=1$ or $5$. Let $C_1$ and $C_2$ be the following collections of $5$-cycles on the set $V$:

(1) Consider a pair of orthogonal $5$-cycle systems of order $10+r$. On $\displaystyle \bigcup_{j=1}^r \{\infty_j\} \cup (\{h_i\}\times \{1,2,3,4,5\})$, when $r=1$ for $i=1,2,...,k$ and when $r=5$ for $i=1$, place a copy of the $5$-cycles of the $5$-cycle system of order $10+r$ in $C_1$ and the $5$-cycles of its orthogonal system in $C_2$.

(2) For the case $r=5$, for each of the remaining holes $h_i, 2 \leq i \leq k$, consider a pair of orthogonal $5$-cycle decompositions of $K_{15}\setminus K_5$. On $\displaystyle \bigcup_{j=1}^5 \{\infty_j\}\cup (h_i\times \{1,2,3,4,5\})$, place a copy of the $5$-cycles in $C_1$ and the $5$-cycles of its orthogonal decomposition in $C_2$, being sure that the vertex set of $K_5$ is $\displaystyle \bigcup_{j=1}^5 \{\infty_j\}$.

(3) If $x$ and $y$ belong to different holes of $H$, place the $5$-cycles $((x,i), (y,i), (x,i+1), (x\circ y,i+3), (y,i+1))\in C_1$ and $((x,i), (y,i), (x,i+2), (x\circ y,i+3), (y,i+2))\in C_2$ for $i=1,2,3,4,5$ and addition is done in modulo $5$.

Then $(V, C_1)$  and $(V, C_2)$ form a pair of orthogonal $5$-cycle systems of order $v$. \hfill $\square$

\section {Orthogonal $6$-cycle systems}

The spectrum of $6$-cycle systems is the set of all $v\equiv 1$ or $9$ (mod $12$). For $v\geq 73$ there exists a pair of orthogonal cyclic $6$-cycle systems of order $12k+1$ \cite{CDDY, CMPP}. The constructions we will give, besides handling the case $v\equiv 9$  (mod $12$), also take care of the missing small systems of the orders $v\equiv 1$  (mod $12$).

\begin{example} \label{ex9}
A pair of orthogonal $6$-cycle systems of order $9$ exists.
\end{example}

Let $V=\mathbb{Z}_9$ and

$C_1= \{(0, 1, 2, 3, 4, 5)$, $(0, 2, 4, 1, 3, 6)$, $(0, 3, 5, 1, 7, 8)$, $(0, 4, 6, 8, 3, 7)$, $(1, 6, 7, 2, 5, 8)$, \linebreak
$(2, 6, 5, 7, 4, 8)\}$,

$C_2= \{(0, 1, 3, 5, 6, 7 )$, $(0, 2, 1, 5, 8, 4)$, $(0, 3, 2, 7, 4, 6)$, $(0, 5, 2, 6, 3, 8)$, $(1, 4, 5, 7, 8, 6)$, \linebreak
$(1, 7, 3, 4, 2, 8)\}$. Then $(V,C_1)$ and $(V,C_2)$ form a pair of orthogonal $6$-cycle systems of order $9$. \hfill $\square$

\begin{example} \label{exK4,4,4}
A pair of orthogonal $6$-cycle decompositions of $K_{4,4,4}$ exists.
\end{example}

Let $V=\{0,1,2,...,11\}$, where $\{0,1,2,3\}$, $\{4,5,6,7\}$, $\{8,9,10,11\}$ be the parts of $K_{4,4,4}$ and

$C_1= \{(0, 4, 1, 5, 2, 6)$, $(0, 5, 3, 4, 2, 7)$, $(0, 8, 1, 6, 3, 9)$, $(0, 10, 1, 7, 3, 11)$, $(1, 9, 5, 10, 4, 11)$, $(2, 8, 7, 11, 6, 9)$, $(2, 10, 6, 8, 5, 11)$, $(3, 8, 4, 9, 7, 10)\}$,

$C_2= \{(0, 4, 2, 8, 1, 11)$, $(0, 5, 1, 9, 2, 10)$, $(0, 6, 8, 3, 5, 9)$, $(0, 7, 3, 10, 5, 8)$, $(1, 4, 3, 6, 9, 7)$, \linebreak
$(1, 6, 2, 11, 7, 10)$, $(2, 5, 11, 4, 8, 7)$, $(3, 9, 4, 10, 6, 11)\}$. Then $(V,C_1)$ and $(V,C_2)$ form a pair of orthogonal $6$-cycle decompositions of $K_{4,4,4}$. \hfill $\square$

\begin{example} \label{ex6-21}
A pair of orthogonal $6$-cycle systems of order $21$ exists.
\end{example}

Let $V= \mathbb{Z}_7 \times \{1,2,3\}$. Let $B_1$ and $B_2$ contain the $5$ base $6$-cycles that will be developed (modulo $7$, --) to produce the $35$ $6$-cycles in $C_1$ and $C_2$, respectively, such that $(V,C_1)$ and $(V,C_2)$ form a pair of orthogonal $6$-cycle systems of order $21$.

$B_1=\{((0,2), (0,1), (1,1), (0,3), (3,1), (6,1))$, $((0,3), (5,2), (2,2), (4,1), (3,2), (1,2))$, \\
$((0,1), (0,3), (4,3), (0,2), (1,3), (2,3))$, $((0,1), (5,1), (1,2), (2,2), (0,3), (5,3))$,\\
$((0,1), (1,3), (5,2), (5,3), (2,1), (4,2))\}$.

$B_2=\{((0,2), (1,1), (3,1), (0,1), (6,1), (0,3))$, $((0,3), (2,2), (3,2), (5,2), (1,2), (4,1))$,\\
$((0,1), (4,3), (1,3), (0,3), (2,3), (0,2))$, $((0,1), (1,2), (0,3), (4,2), (2,1), (2,3))$,\\
$((1,1), (4,2), (1,3), (0,2), (2,1), (0,3))\}$. \hfill $\square$

\begin{lemma} \label{Ort6_1}
A pair of orthogonal $6$-cycle systems of order $v$ exists for all $v\equiv 1, 9, 13$ or $21$ $($mod $24)$.
\end{lemma}

\noindent {\bf Proof} \ A pair of orthogonal $6$-cycle systems of order $13$ is given in \cite{BCP}. The cases $v=9, 21$ are given in Examples \ref{ex9}, \ref{ex6-21} and $v=45$ in the Appendix. The following construction works for all $v=24k+r$ when $r=1, 9$ or $13$ for all $k \geq 1$ and when $r=21$ $k \geq 2$. Let $V=\{\infty\}\cup (X\times \{1,2,3,4\})$, where $|X|=\displaystyle \frac{v-1}{4}$. Let $(X,G,B)$ be a $3-GDD$ of type $2^{(v-1)/8}$ for $r=1, 9$, a $3-GDD$ of type $3^{(v-1)/12}$ for $r=13$ and a $3-GDD$ of type $5^1 3^{2k}$ for $r=21$ \cite{Handbook}. Let $C_1$ and $C_2$ be the following collection of $6$-cycles on
the set $V$:

(1) Consider a pair of orthogonal $6$-cycle systems of order $9, 13$ or $21$ as needed. For each group $g\in G$, on $V=\{\infty\}\cup (g\times \{1,2,3,4\})$ place a copy of the $6$-cycles of the needed system in $C_1$ and the $6$-cycles of its orthogonal system in $C_2$.

(2) For each triple $\{x,y,z\}\in B$, consider a pair of orthogonal $6$-cycle decompositions of $K_{4,4,4}$ (whose parts are $x \times \{1,2,3,4\}$, $y \times \{1,2,3,4\}$, $z \times \{1,2,3,4\}$) on the set $\{x,y,z\} \times \{1,2,3,4\}$, and place a copy of the $6$-cycles of the first system in $C_1$ and the $6$-cycles of its orthogonal system in $C_2$.

Then $(V, C_1)$  and $(V, C_2)$ form a pair of orthogonal $6$-cycle systems of order $v$. \hfill $\square$

\section {Orthogonal $7$-cycle systems}

The spectrum of $7$-cycle systems is the set of all $v\equiv 1$ or $7$ (mod $14$). There exists a pair of orthogonal cyclic $7$-cycle systems of order $56k+15$ for all $k\geq 2$ \cite{CMPP}. Below we will give a construction that works for the orders $14k+1$ and $14k+7$ for all $k\geq 3$.

\begin{example} \label{ex21}
A pair of orthogonal $7$-cycle systems of order $21$ exists.
\end{example}

Let $V=\mathbb{Z}_{20}\cup \{\infty\}$ and $C_1$ contain the $7$-cycles obtained by the two base blocks \\
$(0, 3, 19, 4, 16, 2, 13)$ and $(\infty , 0, 1, 3, 13, 11, 10)$ modulo $20$, where the second one is a short base block.

\noindent Let $C_2 = \{(9, 14, 11, 13, 15, 18, 17)$, $(10, 13, 14, 15, 12, 17, 19)$, $(11, 15, 16, 19, 13, 18, \infty)$, \linebreak
$(0, 8, 4, 11, 1, 10, 12)$, $(0, 11, 2, 10, 3, 7, 15)$, $(0, 14, 1, 9, 2, 12, 16)$, $(0, 17, 1, 12, 3, 8, 18)$, \linebreak
$(0, 19, 1, 13, 2, 14, \infty)$, $(1, 15, 2, 16, 3, 9, 18)$, $(1, 16, 4, 9, 5, 12, \infty)$, $(2, 17, 3, 11, 5, 10, 18)$, \linebreak
$(2, 19, 3, 13, 4, 15, \infty)$, $(3, 14, 4, 12, 6, 7, 18)$, $(3, 15, 5, 13, 6, 8, \infty)$, $(4, 17, 5, 14, 6, 9, 19)$, \linebreak
$(4, 18, 5 , 16, 6, 10, \infty)$, $(5, 19, 6, 11, 7, 9, \infty)$, $(6, 15, 8, 7, 10, 11, 18)$, $(6, 17, 7, 12, 8, 19, \infty)$, \linebreak
$(7, 14, 8, 10, 9, 11, 19)$, $(7, 16, 8, 11, 12, 13, \infty)$, $(8, 13, 9, 15, 10, 16, 17)$, $(9, 12, 14, 10, 17, 11, 16)$, \linebreak
$(0, 1, 2, 3, 4, 5, 6)$, $(0, 2, 4, 1, 3, 5, 7)$, $(0, 3, 6, 1, 7, 4, 10)$, $(0, 4, 6, 2, 5, 8, 9)$, $(0, 5, 1, 8, 2, 7, 13)$, \linebreak
$(12, 18, 16, \infty, 17, 14, 19)$, $(13, 16, 14, 18, 19, 15, 17)\}$. Then $(V,C_1)$ and $(V,C_2)$ form a pair of orthogonal $7$-cycle systems of order $21$. \hfill $\square$

\begin{example} \label{ex21-7}
A pair of orthogonal $7$-cycle decompositions of $K_{21}\setminus K_7$ exists.
\end{example}

Let $V=\mathbb{Z}_{21}$, where the vertex set of the hole $K_7$ is  $V=\{0, 1, 2, 3, 4, 5, 6\}$ and

$C_1= \{(0, 7, 1, 10, 6, 11, 8)$, $(1, 8, 2, 11, 0, 12, 9)$, $(2, 9, 3, 12, 1, 13, 10)$, $(3, 10, 4, 13, 2, 7, 11)$, \linebreak
$(4, 11, 5, 7, 3, 8, 12)$, $(5, 12, 6, 8, 4, 9, 13)$, $(6, 13, 0, 9, 5, 10, 7)$, $(0, 14, 1, 17, 6, 18, 15)$, \linebreak
$(1, 15, 2, 18, 0, 19, 16)$, $(2, 16, 3, 19, 1, 20, 17)$, $(3, 17, 4, 20, 2, 14, 18)$, $(4, 18, 5, 14, 3, 15, 19)$, \linebreak
$(5, 19, 6, 15, 4, 16, 20)$, $(6, 20, 0, 16, 5, 17, 14)$, $(0, 10, 8, 5, 15, 9, 17)$, $(1, 11, 9, 6, 16, 7, 18)$, \linebreak
$(2, 12, 7, 4, 14, 8, 19)$, $(3, 13, 7, 8, 9, 10, 20)$, $(7, 9, 14, 10, 11, 12, 15)$, $(7, 14, 11, 13, 8, 15, 17)$, \linebreak
$(7, 19, 9, 16, 8, 18, 20)$, $(8, 17, 10, 12, 13, 14, 20)$, $(9, 18, 10, 15, 11, 19, 20)$, $(10, 16, 11, 17, 12, 14, 19)$, \linebreak
$(11, 18, 12, 16, 17, 13, 20)$, $(12, 19, 18, 13, 16, 15, 20)$, $(13, 15, 14, 16, 18, 17, 19)\}$,

$C_2= \{(0, 7, 2, 8, 3, 9, 10)$, $(0, 8, 1, 11, 3, 7, 9)$, $(0, 11, 4, 7, 1, 12, 13)$, $(0, 12, 2, 9, 4, 10, 14)$, \linebreak
$(0, 15, 1, 9, 5, 7, 16)$, $(0, 17, 1, 10, 2, 11, 18)$, $(0, 19, 1, 13, 2, 14, 20)$, $(1, 14, 3, 10, 5, 8, 16)$, \linebreak
$(1, 18, 2, 19, 4, 8, 20)$, $(2, 15, 3, 12, 4, 13, 16)$, $(2, 17, 5, 11, 6, 7, 20)$, $(3, 13, 5, 14, 4, 15, 16)$, \linebreak
$(3, 17, 6, 8, 7, 10, 19)$, $(3, 18, 5, 12, 7, 11, 20)$, $(4, 16, 6, 10, 8, 9, 20)$, $(4, 17, 7, 13, 6, 9, 18)$, \linebreak
$(5, 15, 7, 14, 6, 12, 19)$, $(5, 16, 9, 11, 8, 12, 20)$, $(6, 15, 8, 14, 9, 12, 18)$, $(6, 19, 7, 18, 10, 13, 20)$, \linebreak
$(8, 13, 9, 15, 10, 16, 18)$, $(8, 17, 11, 13, 15, 18, 19)$, $(9, 17, 10, 11, 14, 16, 19)$, $(10, 12, 14, 18, 17, 16, 20)$, \linebreak
$(11, 12, 17, 14, 13, 19, 15)$, $(11, 16, 12, 15, 20, 17, 19)$, $(13, 17, 15, 14, 19, 20, 18)\}$. Then $(V,C_1)$ and $(V,C_2)$ form a pair of orthogonal $7$-cycle decompositions of $K_{21}\setminus K_7$. \hfill $\square$

\begin{lemma} \label{Ort7}
A pair of orthogonal $7$-cycle systems of order $v$ exists for all $v \equiv 1$ or $7$ $($mod $14)$, except $v=7$.
\end{lemma}

\noindent {\bf Proof} \ A pair of orthogonal $7$-cycle systems of order $7$ does not exist since a $7$-cycle system of order $7$ has $3$ blocks. A pair of orthogonal $7$-cycle systems of order $15$ is given in \cite{BCP}. The case $v=21$ is given in Example \ref{ex21} and the cases $v=29,35$ are given in the Appendix. The following construction works for all $v=14k+r$ and all $k \geq 3$, where $r=1$ or $7$. Let $(Q, \circ )$ be a commutative quasigroup of order $2k$, for $k\geq 3$, with holes $H=\{h_1,h_2,...,h_k\}$, where $h_i=\{2i-1,2i\}$ for $i=1,2,...,k$. Let $V=\displaystyle \bigcup_{j=1}^r \{\infty_j\} \cup (Q\times \{1,2,3,4,5,6,7\})$, where $r=1$ or $7$. Let $C_1$ and $C_2$ be the following collections of $7$-cycles on the set $V$:

(1) Consider a pair of orthogonal $7$-cycle systems of order $14+r$. On $\displaystyle \bigcup_{j=1}^r \{\infty_j\} \cup (\{h_i\}\times \{1,2,3,4,5,6,7\})$ when $r=1$ for $i=1,2,...,k$ and when $r=7$ for $i=1$, place a copy of the $7$-cycles of the system of order $14+r$ in $C_1$ and the $7$-cycles of its orthogonal system in $C_2$.

(2) When $r=7$, for each of the remaining holes $h_i$ with $i=2,3,...,k$, consider a pair of orthogonal $7$-cycle decompositions of $K_{21}\setminus K_7$. On $\displaystyle \bigcup_{j=1}^7 \{\infty_j\} \cup (\{h_i\}\times \{1,2,3,4,5,6,7\})$ for $i=2,3,...,k$, place a copy of the $7$-cycles of the first system in $C_1$ and the $7$-cycles of its orthogonal system in $C_2$, being sure that the vertex set of $K_7$ is $\displaystyle \bigcup_{j=1}^7 \{\infty_j\}$.

(3) If $x$ and $y$ belong to different holes of $H$ with $\{x,x'\}\in H$ and $\{y,y'\}\in H$. Then place the $7$-cycles $((x,i), (y,i), (x,i+1), (y,i+3), (x\circ y,i+6), (x,i+3), (y,i+1))\in C_1$ and $((x,i), (y,i), (x',i+3), (y,i+4), (x\circ y,i+6), (x,i+4), (y',i+3))\in C_2$ for $i=1,2,3,4,5,6,7$ modulo $7$.

Then $(V, C_1)$  and $(V, C_2)$ form a pair of orthogonal $7$-cycle systems of order $v$. \hfill $\square$

\section {Orthogonal $8$-cycle systems}

The spectrum of $8$-cycle systems is the set of all $v\equiv 1$ (mod $16$). There exists a pair of orthogonal $8$-cycle systems of order $16k+1$ for all $16k+1\geq 129$ \cite{BCP}.

The construction we will give to obtain the missing small orders work also for all $v=16k+1 \geq 33$.

\begin{example} \label{exK16,16}
A pair of orthogonal $8$-cycle decompositions of $K_{16,16}$ exists.
\end{example}
Let $X=\{(i,0): i\in \mathbb{Z}_{16}\}$ and $Y=\{(i,1): i\in \mathbb{Z}_{16}\}$ be the parts of $K_{16,16}$, $C_1$ and $C_2$, each, be a collection of $8$-cycles obtained by developing the base $8$-cycles $((0,0),(0,1),(1,0),(2,1),(6,0)$, \\
$(1,1),(4,0),(6,1))$ and $((0,0),(0,1),(6,0),(7,1),(4,0),(3,1),(7,0),(5,1))$, respectively (mod $16$, --).  Then $(V,C_1)$ and $(V,C_2)$ form a pair of orthogonal $8$-cycle decompositions of $K_{16,16}$. \hfill $\square$

\begin{lemma} \label{Ort8}
A pair of orthogonal $8$-cycle systems of order $v$ exists for all $v \equiv 1$ $($mod $16)$.
\end{lemma}

\noindent {\bf Proof} \ A pair of orthogonal $8$-cycle systems of order $17$ is given in \cite{BCP}. For all $v=16k+1 \geq 33$, let $X=\{\infty \} \cup \{(i,j)\ | \ 1\leq i\leq k, 1\leq j\leq 16\}$.

\noindent (1) Consider a pair of orthogonal $8$-cycle systems of order $17$. For $1\leq i\leq k$, on each set $\{\infty \} \cup \{(i,j)\ | \
1\leq j\leq 16\}$, place a copy of the $8$-cycles of the first system in $C_1$  and the $8$-cycles of its orthogonal system in $C_2$.

\noindent (2) For each pair of $x, y \in \{1,2,...,k\}$, consider a pair of orthogonal $8$-cycle decompositions of $K_{16,16}$ on $\{(x,j)\ | \ 1\leq j\leq 16\} \cup \{(y,j)\ | \ 1\leq j\leq 16\}$ and place the $8$-cycles of the first system in $C_1$ and the $8$-cycles of its orthogonal system in  $C_2$.

Then $(V, C_1)$  and $(V, C_2)$ form a pair of orthogonal $8$-cycle systems of order $v$. \hfill $\square$

\section {Orthogonal $9$-cycle systems}

The spectrum of $9$-cycle systems is the set of all $v\equiv 1$ or $9$ (mod $18$). There exists a pair of orthogonal cyclic $9$-cycle systems of order $72k+55$ for all $72k+55\geq 199$ \cite{CMPP}. The constructions we will give work for all $v\equiv 1$ or $9$  (mod $18$) when $v\geq 55$.

\begin{example} \label{ex27}
A pair of orthogonal $9$-cycle systems of order $27$ exists.
\end{example}

Let $V=\mathbb{Z}_{26}\cup \{\infty\}$ and

$C_1=\{(i, 1+i, 3+i, 6+i, 10+i, 2+i, 18+i, 24+i, 5+i): 0\leq i \leq 26 \}\cup \{(\infty , 2j, 9+2j, 20+2j, 6+2j, 19+2j, 2+2j, 13+2j, 25+2j): 0\leq j \leq 13  \} $ modulo $26$,

$C_2=\{(i, 5+i, 3+i, 9+i, 25+i, 6+i, 14+i, 2+i, 11+i): 0\leq i \leq 26 \}\cup \{(\infty , 2j, 1+2j, 5+2j, 8+2j, 21+2j, 18+2j, 14+2j, 13+2j): 0\leq j \leq 13 \} $ modulo $26$. Then $(V,C_1)$ and $(V,C_2)$ form a pair of orthogonal $9$-cycle systems of order $27$. \hfill $\square$

\begin{example} \label{ex27-9}
A pair of orthogonal $9$-cycle decompositions of $K_{27}\setminus K_9$ exists.
\end{example}

Let $V=\mathbb{Z}_{27}$, where the vertex set of the hole $K_9$ is  $V=\{0, 1, 2, 3, 4, 5, 6, 7, 8\}$ and

$C_1= \{(0, 9, 1, 12, 6, 13, 3, 15, 14)$, $(1, 10, 2, 13, 7, 14, 4, 16, 15)$, $(2, 11, 3, 14, 8, 15, 5, 17, 16)$,\linebreak
$(3, 12, 4, 15, 0, 16, 6, 9, 17)$, $(4, 13, 5, 16, 1, 17, 7, 10, 9)$, $(5, 14, 6, 17, 2, 9, 8, 11, 10)$, \linebreak
$(6, 15, 7, 9, 3, 10, 0, 12, 11)$, $(7, 16, 8, 10, 4, 11, 1, 13, 12)$, $(8, 17, 0, 11, 5, 12, 2, 14, 13)$, \linebreak
$(0, 18, 1, 21, 6, 22, 3, 24, 23)$, $(1, 19, 2, 22, 7, 23, 4, 25, 24)$, $(2, 20, 3, 23, 8, 24, 5, 26, 25)$, \linebreak
$(3, 21, 4, 24, 0, 25, 6, 18, 26)$, $(4, 22, 5, 25, 1, 26, 7, 19, 18)$, $(5, 23, 6, 26, 2, 18, 8, 20, 19)$, \linebreak
$(6, 24, 7, 18, 3, 19, 0, 21, 20)$, $(7, 25, 8, 19, 4, 20, 1, 22, 21)$, $(8, 26, 0, 20, 5, 21, 2, 23, 22)$, \linebreak
$(9, 19, 12, 23, 17, 21, 24, 20, 18)$, $(10, 20, 13, 24, 9, 22, 25, 21, 19)$, $(11, 21, 14, 25, 10, 23, 26, 22, 20)$, \linebreak
$(12, 22, 15, 26, 11, 24, 18, 23, 21)$, $(13, 23, 16, 18, 12, 25, 19, 24, 22)$, $(14, 24, 17, 19, 13, 26, 20, 25, 23)$, \linebreak
$(15, 25, 9, 20, 14, 18, 21, 26, 24)$, $(16, 26, 10, 21, 15, 19, 22, 18, 25)$, $(17, 18, 11, 22, 16, 20, 23, 19 26)$, \linebreak
$(0, 13, 9, 5, 18, 10, 12, 17, 22)$, $(1, 14, 9, 11, 7, 20, 12, 15, 23)$, $(2, 15, 9, 12, 8, 21, 13, 10, 24)$, \linebreak
$(3, 16, 9, 23, 11, 13, 15, 17, 25)$, $(4, 17, 10, 6, 19, 14, 16, 12, 26)$, $(9, 21, 16, 11, 15, 10, 22, 14, 26)$, \linebreak
$(10, 14, 11, 17, 20, 15, 18, 13, 16)$, $(11, 19, 16, 24, 12, 14, 17, 13, 25)\}$,

$C_2= \{(0, 9, 19, 2, 13, 5, 21, 20, 23)$, $(1, 10, 20, 3, 14, 6, 22, 21, 24)$, $(2, 11, 21, 4, 15, 7, 23, 22, 25)$, \linebreak
$(3, 12, 22, 5, 16, 8, 24, 23, 26)$, $(4, 13, 23, 6, 17, 0, 25, 24, 18)$, $(5, 14, 24, 7, 9, 1, 26, 25, 19)$, \linebreak
$(6, 15, 25, 8, 10, 2, 18, 26, 20)$, $(7, 16, 26, 0, 11, 3, 19, 18, 21)$, $(8, 17, 18, 1, 12, 4, 20, 19, 22)$, \linebreak
$(0, 14, 20, 25, 16, 9, 17, 1, 13)$, $(1, 15, 21, 26, 17, 10, 9, 2, 14)$, $(2, 16, 22, 18, 9, 11, 10, 3, 15)$, \linebreak
$(3, 17, 23, 19, 10, 12, 11, 4, 16)$, $(4, 9, 24, 20, 11, 13, 12, 5, 17)$, $(5, 10, 25, 21, 12, 14, 13, 6, 9)$, \linebreak
$(6, 11, 26, 22, 13, 15, 14, 7, 10)$, $(0, 21, 2, 22, 11, 1, 25, 12, 15)$, $(1, 22, 3, 23, 12, 2, 26, 13, 16)$, \linebreak
$(2, 23, 4, 24, 13, 3, 18, 14, 17)$, $(3, 24, 5, 25, 14, 4, 19, 15, 9)$, $(4, 25, 6, 26, 15, 5, 20, 16, 10)$, \linebreak
$(5, 26, 7, 18, 16, 6, 21, 17, 11)$, $(6, 18, 8, 19, 17, 7, 22, 9, 12)$, $(7, 19, 0, 20, 9, 8, 23, 10, 13)$, \linebreak
$(8, 20, 1, 21, 10, 0, 24, 11, 14)$, $(0, 12, 7, 11, 8, 13, 9, 23, 16)$, $(0, 18, 11, 15, 16, 19, 12, 20, 22)$, \linebreak
$(1, 19, 6, 24, 16, 21, 3, 25, 23)$, $(2, 20, 7, 25, 9, 26, 4, 22, 24)$, $(5, 18, 13, 17, 16, 14, 9, 21, 23)$, \linebreak
$(8, 12, 18, 15, 23, 11, 19, 14, 26)$, $(8, 15, 22, 10, 26, 12, 24, 19, 21)$, $(10, 14, 22, 17, 25, 13, 19, 26, 24)$, \linebreak
$(10, 15, 17, 20, 13, 21, 14, 23, 18)$, $(11, 16, 12, 17, 24, 15, 20, 18, 25)\}$. Then $(V,C_1)$ and $(V,C_2)$ form a pair of orthogonal $9$-cycle decompositions of $K_{27}\setminus K_9$. \hfill $\square$

\begin{example} \label{ex37}
A pair of orthogonal $9$-cycle systems of order $37$ exists.
\end{example}

Let $V= \mathbb{Z}_{37}$, $B_1$ and $B_2$ contain the two $9$-cycles that will be developed modulo $37$ to produce the $74$ \ $9$-cycles in $C_1$ and $C_2$.

$B_1=\{(0, 2, 5, 9, 14, 20, 3, 10, 29)$, $(0, 9, 19, 30, 5, 18, 32, 17, 1)\}$.

$B_2=\{(0, 2, 11, 14, 33, 6, 10, 24, 30)$, $(0, 5, 16, 29, 14, 26, 10, 30, 1)\}$. Then $(V,C_1)$ and $(V,C_2)$ form a pair of orthogonal $9$-cycle systems of order $37$. \hfill $\square$

\begin{example} \label{ex9-45}
A pair of orthogonal $9$-cycle systems of order $45$ exists.
\end{example}

Let $V= \mathbb{Z}_{11} \times \{1,2,3,4\}\cup \{\infty\}$. Let $B_1$ and $B_2$ contain the $10$ base $9$-cycles that will be developed modulo ($11$, --) to produce the $110$ $9$-cycles in $C_1$ and $C_2$.

$B_1=\{(\infty, (1,2), (8, 2), (0,1), (7,3), (10, 1), (0, 2), (5,4), (8,4))$, $(\infty,(8,1),(3,2),(0,3),(10,2)$,\linebreak
$(6,4),(0,4),(4,1),(5,3))$, $((0,1), (1, 1), (1, 2), (2, 1), (2,4), (0,4), (5, 1), (8,3), (5, 3))$, $((0,1),(4,1)$, \linebreak
$(2,4),(1,1),(3,2),(5,2),(5,4),(1,4),(10,3))$, $((0,2),(7,1),(10,1),(6,2),(1,3),(8,3),(7,4)$, \linebreak
$(10,3), (10,4))$, $((0,3), (10,3), (10, 1), (1, 3), (1, 2), (6, 2), (4,3), (9,1), (1,4))$, $((0,4),(4,3),(0,1)$, \linebreak
$(2,4), (10,3),(8,2),(6,4),(2,1),(10,4))$, $((0,4),(7,3),(4,2),(1,4),(8,2),(10,1),(5,1),(7,1)$, \linebreak
$(5,3))$, $((0,2),(6,4),(8,3),(1,2),(3,4),(9,1),(3,2),(7,3),(1,4))$, $((0,2),(8,1),(7,4),(4,2)$, \linebreak
$(1,2),(2,2),(7,3),(1,3),(10,3))\}$.

$B_2=\{((1,2),(8,2),(5,3),(8,3),(4,4),(5,2),(10,2),(8,4),(4,2))$, $((3,2),(0,1),(8,2),(2,1)$, \linebreak
$(4,2),(10,1),(5,1),(3,1),(5,3))$, $((0,1),(7,3),(6,1),(9,3),(9,1),(2,3),(7,4),(8,1),(6,3))$, \linebreak
$(10,1),(7,3),(2,1),(1,3),(0,2),(6,3),(8,2),(8,3),(6,4))$, $((10,1),(0,2),(7,4),(1,1),(10,4)$, \linebreak
$(2,3),(3,3),(9,4),(8,2))$, $((0,2),(5,4),(10,4),(1,4),(10,3),(9,4),(8,3),(0,4),(9,2))$, $((0,1)$, \linebreak
$(1,1),(8,2),(0,3),(5,3),(1,2),(1,4),(4,4),(4,1))$, $(\infty,(1,1),(2,4),(10,1),(7,4),(5,1),(8,1)$, \linebreak
$(8,2),(9,2))$, $((1,4),(0,4),(7,1),(0,2),(3,4),(6,2),(2,3),(4,3),(8,3))$, $((0,1),(5,4),(9,4),\infty$, \linebreak
$(0,3),(0,4),(5,2),(4,3),(10,2))\}$. Then $(V,C_1)$ and $(V,C_2)$ form a pair of orthogonal $9$-cycle systems of order $45$. \hfill $\square$

\newpage

\begin{example} \label{exK9,9,9}
A pair of orthogonal $9$-cycle decompositions of $K_{9,9,9}$ exists.
\end{example}

Let $i,j$ be integers and $\{(i,j): 0\leq i\leq 8\}$ be the groups of $K_{9,9,9}$  for $0\leq j\leq 2$. Let \\
$C_1=\{((0+i,0+j),(0+i,1+j),(1+i,2+j),(3+i,0+j),(6+i,1+j),(2+i,2+j),(6+i,0+j),(4+i,1+j),(3+i,2+j)): 0\leq i\leq 8  \ \text{and} \  0\leq j\leq 2 \}$ and \\ $C_2=\{((0+i,0+j),(0+i,1+j),(2+i,2+j),(7+i,0+j),(5+i,1+j),(6+i,2+j),(3+i,0+j),(7+i,1+j),(1+i,2+j)): 0\leq i\leq 8  \ \text{and} \  0\leq j\leq 2 \}$, (mod $9$, mod $3$).  Then $(V,C_1)$ and $(V,C_2)$ form a pair of orthogonal $9$-cycle decompositions of $K_{9,9,9}$. \hfill $\square$

\begin{lemma} \label{Ort9}
A pair of orthogonal $9$-cycle systems of order $v$ exists $v \equiv 1$ or $9$ $($mod $18)$, except $v=9$.
\end{lemma}

\noindent {\bf Proof} \ A pair of orthogonal $9$-cycle systems of order $9$ does not exist since a $9$-cycle system of order $9$ has $4$ blocks. A pair of orthogonal $9$-cycle systems of order $19$ is given in \cite{BCP}. The case $v=27,37$ and $45$ are given in Examples \ref{ex27}, \ref{ex37} \ref{ex9-45}. For all $v=18k+r$ and all $k \geq 3$, where $r=1$ or $9$, let $V=\displaystyle \bigcup_{j=1}^r \{\infty_j\} \cup (X\times \{l: 1\leq l \leq 9\})$, where $|X|=2k$. Let $(X,G,B)$ be a $3-GDD$ of type $2^k$ when $2k\equiv 0$ or $2$ (mod $6$) and a $3-GDD$ of type $4^12^{k-2}$ when $2k\equiv 4$ (mod $6$). They exist for $k\geq 3$ \cite{Handbook}. Let $C_1$ and $C_2$ be the following collection of $9$-cycles on
the set $V$:

(1) Consider a pair of orthogonal $9$-cycle systems of order $19, 27, 37$ and $45$ as needed.

On $\displaystyle \bigcup_{j=1}^r \{\infty_j\} \cup (g_1\times \{l: 1\leq l \leq 9\})$, place a copy of the $9$-cycles of the needed system in $C_1$ and the $9$-cycles of its orthogonal system in $C_2$.

(2) Consider a pair of orthogonal $9$-cycle decompositions of $K_{27}\setminus K_r$. For each group $g_i\in G$, when  $i\geq 2$, on $\displaystyle \bigcup_{j=1}^r \{\infty_j\} \cup (\{g_i\}\times \{l: 1\leq l \leq 9\})$, place a copy of the $9$-cycles of the needed decomposition in $C_1$ and the $9$-cycles of its orthogonal decomposition in $C_2$, being sure that the vertex set of $K_9$ is $\displaystyle \bigcup_{j=1}^9 \{\infty_j\}$.

(3) For each triple $\{x,y,z\}\in B$, on $\{x,y,z\} \times \{l: 1\leq l \leq 9\}$, consider a pair of orthogonal $9$-cycle decompositions of $K_{9,9,9}$ whose parts are $x \times \{l: 1\leq l \leq 9\}$, $y \times \{l: 1\leq l \leq 9\}$, $z \times \{l: 1\leq l \leq 9\}$, and place a copy of the $9$-cycles in $C_1$ and the $9$-cycles of its orthogonal system in $C_2$.

Then $(V, C_1)$  and $(V, C_2)$ form a pair of orthogonal $9$-cycle systems of order $v$. \hfill $\square$

Combining Lemmas \ref{Ort5}, \ref{Ort6_1}, \ref{Ort7}, \ref{Ort8} and \ref{Ort9} gives the following result.

\begin{theorem}
A pair of orthogonal $l$-cycle systems of order $v$ exists for all admissible $v$ when $l=5,6,7, 8\ or\ 9$, except $7$-cycle system of order $7$ and $9$-cycle system of order $9$, except $(l,v)=(7,7)$ and $(l,v)=(9,9)$.
\end{theorem}

\newpage

\section {Appendix}

In the following cases $(V,C_1)$ and $(V,C_2)$ will form a pair of orthogonal $l$-cycle systems of order $v$ for considered $l$ and $v$.

\begin{example} \label{ex5-21}
A pair of orthogonal $5$-cycle systems of order $21$ exists.
\end{example}

Let $V=\mathbb{Z}_{21}$, $B_1$ and $B_2$ contain the two $5$-cycles that will be developed modulo $21$ to produce the $5$-cycles in $C_1$ and $C_2$.

$B_1=\{(0, 1, 3, 6, 12)$, $(0, 4, 9, 17, 10)\}$.

$B_2=\{(0, 1, 4, 13, 15)$, $(0, 4, 12, 2, 7)\}$. \hfill $\square$

\begin{example} \label{ex25}
A pair of orthogonal $5$-cycle systems of order $25$ exists.
\end{example}

Let $V= \mathbb{Z}_5 \times \mathbb{Z}_5$ and $B_1$ and $B_2$ each contains the $4$ base $5$-cycles that the first two base cycles will be developed modulo $(5,5)$ and the last two base cycles will be developed modulo $(5, -)$ to produce the $60$ $5$-cycles in each of $C_1$ and $C_2$.

$B_1=\{((0,0), (1,0), (4,1), (3,3), (1,1))$, $((0,0), (2,0), (2,1), (2,3), (4,1))$, \\
$((0,0), (1,2), (2,4), (3,1), (4,3))$, $((0,0), (2,1), (4,2), (1,3), (3,4))\}$.

$B_2=\{((0,0), (1,0), (4,2), (3,3), (1,2))$, $((0,0), (2,0), (2,2), (2,3), (4,2))$, \\
$((0,0), (1,1), (2,2), (3,3), (4,4))$, $((0,0), (2,2), (4,4), (1,1), (3,3))\}$. \hfill $\square$

\begin{example} \label{exK6,10}
A pair of orthogonal $6$-cycle decompositions of $K_{6,10}$ exists.
\end{example}

Let $X=\{0,1,2,3,4,5\}$ and $Y=\{6,7,...,15\}$ be the parts of $K_{6,10}$ and

$C_1= \{(0, 6, 1, 7, 2, 8)$, $(0, 7, 3, 6, 2, 9)$, $(0, 10, 1, 8, 3, 11)$, $(0, 12, 1, 9, 3, 13)$, $(0, 14, 1, 11, 2, 15)$, $(1, 13, 5, 7, 4, 15)$, $(2, 10, 4, 9, 5, 12)$, $(2, 13, 4, 6, 5, 14)$, $(3, 10, 5, 11, 4, 12)$, $(3, 14, 4, 8, 5, 15)\}$,

$C_2= \{(0, 6, 2, 10, 1, 13)$, $(0, 7, 1, 8, 4, 12)$, $(0, 8, 3, 6, 4, 14)$, $(0, 9, 1, 6, 5, 10)$, $(0, 11, 4, 7, 3, 15)$, $(1, 11, 3, 9, 5, 15)$, $(1, 12, 2, 7, 5, 14)$, $(2, 8, 5, 12, 3, 13)$, $(2, 9, 4, 13, 5, 11)$, $(2, 14, 3, 10, 4, 15)\}$. \hfill $\square$

\newpage

\begin{example} \label{ex45}
A pair of orthogonal $6$-cycle systems of order $45$ exists.
\end{example}

Let $V=\{\infty\} \cup \displaystyle \bigcup_{i=0}^3 X_i \cup \displaystyle \bigcup_{j=0}^1 Y_j$, where $X_i=\{x_{6i+1}, x_{6i+2}, ..., x_{6i+6}\}$ and $Y_j=\{y_{10j+1}$, $y_{10j+2}$, ..., $y_{10j+10}\}$. Place a pair of orthogonal $6$-cycle systems of order $25$ and $21$ on $\{\infty\} \cup \displaystyle \bigcup_{i=0}^3 X_i$ and $\{\infty\} \cup  \displaystyle \bigcup_{j=0}^1 Y_j$, respectively.

Finally, for each $i=0,1,2,3$, $j=0,1$, place a pair of orthogonal $6$-cycle decompositions of $K_{6,10}$ on $X_i \cup Y_j$. \hfill $\square$

\begin{example} \label{ex29}
A pair of orthogonal $7$-cycle systems of order $29$ exists.
\end{example}

Let $V= \mathbb{Z}_{29}$, $B_1$ and $B_2$ contain the two $7$-cycles that will be developed modulo $29$ to produce the $58$ \ $7$-cycles in $C_1$ and $C_2$.

$B_1=\{(0, 1, 3, 6, 10, 15, 21)$, $(0, 7, 16, 26, 15, 27, 14)\}$.

$B_2=\{(0, 4, 12, 21, 7, 10, 17)$, $(0, 1, 6, 16, 3, 5, 11)\}$. \hfill $\square$

\begin{example} \label{ex35}
A pair of orthogonal $7$-cycle systems of order $35$ exists.
\end{example}

Let $V=\mathbb{Z}_{34}\cup \{\infty\}$ and

$C_1=\{(i, 3+i, 10+i, 21+i, 19+i, 27+i, 14+i), (i, 1+i, 5+i, 10+i, 19+i, 31+i, 16+i): 0\leq i \leq 34 \}\cup \{(\infty , 2j, 6+2j, 16+2j, 33+2j, 5+2j, 15+2j): 0\leq j \leq 17  \} $ modulo $34$,

$C_2=\{(i, 7+i, 19+i, 30+i, 22+i, 28+i, 15+i), (i, 4+i, 6+i, 15+i, 20+i, 30+i, 16+i): 0\leq i \leq 34 \}\cup \{(\infty , 2j, 1+2j, 4+2j, 21+2j, 18+2j, 17+2j): 0\leq j \leq 17 \} $ modulo $34$. \hfill $\square$

\end{document}